\documentclass[aps,eqsecnum,nofootinbib,superscriptaddress]{revtex4}
\usepackage[dvips]{graphicx}  % for graphics
\usepackage{amsmath,amsfonts}
\usepackage{color}  % for color
\def\be{\begin{eqnarray}}
\def\ee{\end{eqnarray}}
\def\beq{\begin{equation}}
\def\eeq{\end{equation}}

\def\({\left (}
\def\){\right )}

\newtheorem{theorem}{Theorem}[section]
\newtheorem{lemma}[theorem]{Lemma}

\newtheorem{definition}{Definition}[section]

\newcommand{\qed}{\nobreak \ifvmode \relax \else
      \ifdim\lastskip<1.5em \hskip-\lastskip
      \hskip1.5em plus0em minus0.5em \fi \nobreak
      \vrule height0.75em width0.5em depth0.25em\fi}
\newcommand{\qcd}{\begin{flushright} $\Box$ \end{flushright}}
\begin{document}
%%%%%%%%%%%%%%%%%%%%%%%%%%%%%%%%%%%%%%%%%%%%%%%%%%%%%%%%%%%%%
%%%%%%%%%%%%%%%%%%%%%%   Title Page   %%%%%%%%%%%%%%%%%%%%%%%
%\rightline{Preprint-MTM}
%\rightline{hep-th/0605224}

\title{Compact Mean Convex Hypersurfaces and the Fundamental Group of Manifolds with Nonnegative Ricci Curvature}

\author{I.P. Costa e Silva}
\email{ivanpcs@mtm.ufsc.br}
\affiliation{Department of Mathematics,\\ 
Universidade Federal de Santa Catarina \\88.040-900 Florian\'{o}polis-SC, Brasil}

\date{\today}

\begin{abstract}
We show that the existence of an embedded compact, boundaryless hypersurface $S$ of strictly positive mean curvature in a noncompact, connected, complete Riemannian $n$-manifold $N$ of nonnegative Ricci curvature ($n \geq 2$) implies that the map $i_{\#}:\pi_1(S) \rightarrow \pi_1(N)$ induced by the inclusion is surjective, provided only that $N\setminus S$ has two connected components $N_{\pm}$ with $S\cup N_{+}$ noncompact and $i_{\#}:\pi_1(S) \rightarrow \pi_1(S\cup N_{+})$ is surjective. The idea of the proof is simple: to view $N$ as a spacelike hypersurface in a suitable Lorentz manifold and then apply a recent version \cite{GLme} of certain classic results by Gannon \cite{gannon1,gannon2} and Lee \cite{lee} on the topology of spacetimes. As an application, we show that if $N$ has an asymptotically flat end then this is its {\em only} end, and $N$ is simply connected for $n \geq 3$.  
\end{abstract}
%
%\pacs{04.20.Dw;04.20.Gz;02.40.-k}
\maketitle

%%%%%%%%%%%%%%%%%%%%%%%%%%%%%%%%%%%%%%%%%%%%%%%%%%%%%%%%%%%%%
\section{Introduction}\label{sec:intro}
%%%%%%%%%%%%%%%%%%%%%%%%%%%%%%%%%%%%%%%%%%%%%%%%%%%%%%%%%%%%%

Let $N$ be a complete, noncompact, connected $n$-dimensional ($n \geq 2$) Riemannian manifold with Ricci curvature $Ric_N\geq 0$. As early as 1968, J. Milnor \cite{milnor} conjectured that the fundamental group of $N$ must be finitely generated. (If $N$ is compact, then Morse theory \cite{morse} guarantees that $N$ is homotopically a $CW$-complex with finitely many cells in any dimension, and therefore $\pi_1(N)$ is indeed finitely generated in this case.) Although it is still an open problem whether this conjecture holds in general, much progress has been made to establish it in special cases, in which a large number of different techniques are brought to bear (see, for example, \cite{wei,sormanisurvey} for recent surveys on manifolds of nonnegative Ricci curvature in general and on Milnor's conjecture in particular, with extensive references). Perhaps the most well-known case is when there are bounds on the diameter or volume growth of $N$ from comparison with Euclidean space. Let $p \in N$ and let $B_p(r)$ be the geodesic open ball of radius $r>0$ centered at $p$. Denote the Euclidean volume of the unit ball in $\mathbb{R}^n$ by $\omega_n$. Li \cite{li} and Anderson \cite{anderson} have independently shown that if the volume growth 
\[
\alpha_N := \lim \inf_{r \rightarrow \infty} \frac{\mbox{vol}(B_p(r))}{\omega _nr^n}
\]
of $N$ (compared to Euclidean space $E^n=(\mathbb{R}^n,g_E)$) is Euclidean (or large), i.e., $\alpha_N >0$, then the fundamental group is finitely generated. Improving on this result, in 1994 Perelman \cite{perelman} proved that there exists a constant $0<\delta_n<1$ depending only on the dimension $n$ such that if $\alpha_N \geq 1-\delta_n$ then $N$ is contractible. Cheeger and Colding \cite{CC} later established  that for a (presumably larger) constant for which $0<\delta'_n<1$ and $\alpha_N \geq 1-\delta'_n$, $N$ is diffeomorphic to $\mathbb{R}^n$. As for diameter bounds, Sormani \cite{sormani} has shown that if $N$ has small linear diameter growth:
\[
\lim \sup_{r \rightarrow \infty} \frac{\mbox{diam}(\partial B_p(r))}{\omega _nr^n} <s_n
\]
for an explicitly given $s_n$, then $\pi_1(N)$ is finitely generated. 

The purpose of this note is to use an entirely different approach to the general problem of relating $\pi_1(N)$ and the geometry of $N$. It differs from the previous geometric comparison techniques in at least two different ways. First, a comparison with $E^n$ is made in an indirect manner, by showing that if $N$ admits an embedded codimension one submanifold with certain features which resemble the standard embedding of $S^{n-1}$ in $\mathbb{R}^n$, then this will place constraints on the topology of $N$. Specifically our main theorem is the following:  

\begin{theorem}
\label{maintheorem}
Let $N$ be a complete, connected $n$-dimensional ($n \geq 2$) Riemannian manifold with Ricci curvature $Ric_N\geq 0$. Let $S \subset N$ be a smooth codimension one embedded submanifold such that 
\begin{itemize}
\item[(1)] $S$ is compact without boundary, connected and has positive mean curvature $H_S >0$ [with respect to a suitable normal - see conventions below].
\item[(2)] $N\setminus S$ has two connected components $N_{\pm}$ with $S \cup N_{+}$ noncompact and $i_{\#}:\pi_1(S) \rightarrow \pi_1(S \cup N_{+})$ induced by the inclusion $i:S \hookrightarrow N$ is surjective. 
\end{itemize}
Then $S \cup N_{-}$ is compact and $i_{\#}:\pi_1(S) \rightarrow \pi_1(N)$ is surjective. In particular, $\pi_1(N)$ is finitely generated. 
\end{theorem}    

Secondly, in the proof of this result we view $N$ as a {\em spacelike submanifold} of a suitable Lorentzian manifold, and use certain known theorems about these manifolds. The techniques used in proving the latter theorems include only standard methods in causal theory, and relatively simple results about covering manifolds. Thus our proof is entirely geometric, whereas the comparison results mentioned above often require delicate analytical arguments. 

The price of the above-mentioned simplicity is that our results are restricted to $\pi_1(N)$. This can give good topological information in dimension $n=3$, but in higher dimensions it is (in general) a somewhat poor invariant. The topological structure of $N$ if it has large volume growth but violates the Perelman bound can be quite involved even with small diameter growth. For instance, a famous example constructed by Menguy for $n=4$ \cite{menguy} has infinite topological type. But since $\pi_1(N)$ is easier to control, results about its structure are frequently easier to obtain.  

The rest of the paper is organized as follows. We first give some preliminary definitions in Section \ref{preliminaries} to set the conventions, general assumptions and notation, and some results on Lorentzian manifolds pertinent to our proof. Section \ref{AsympFlat} gives an important application of our main result: we show that if (together with the above assumptions) $N$ has an {\em asymptotically flat end} and $n \geq 3$, then this is its {\em only} end, and it is simply connected (in the case $n=3$ it is diffeomorphic to $\mathbb{R}^3$). We end the Section with a few relevant remarks on the latter result. The proof of the main result is deferred to Section \ref{detailedproofs}. 

%%%%%%%%%%%%%%%%%%%%%%%%%%%%%%%%%%%%%%%%%%%%%%%%%%%%%%%%%%%%%
\section{Preliminaries: Basic Definitions \& Results about Spacetimes}\label{preliminaries}
%%%%%%%%%%%%%%%%%%%%%%%%%%%%%%%%%%%%%%%%%%%%%%%%%%%%%%%%%%%%%

 In this Section, we fix an $(n+1)$-dimensional connected manifold $M$ endowed with a smooth \footnote{In this paper we always use `smooth' to mean $C^{\infty}$. Moreover, 'manifold' always means a Hausdorff, second countable $C^{\infty}$ manifold without boundary, but no {\em a priori} assumptions about orientability are made.} Lorentzian metric tensor $g$ (signature $(-,+, \ldots,+)$). We briefly review here some pretty standard definitions and certain pertinent results in causal theory of spacetimes, mostly to set the terminology; more details can be found in the core references \cite{HE,wald,oneill,BE}. 

A nonzero vector $v \in TM$ is said to be {\em timelike} [resp. {\em lightlike} (or {\em null}), {\em spacelike}, {\em nonspacelike} (or {\em causal})] if $g(v,v)<0$ [resp. $g(v,v) =0$, $g(v,v) >0$, $g(v,v) \leq 0$]. A differentiable curve $\gamma:I \subseteq \mathbb{R}\rightarrow M$ (where $I$ denotes some interval) is said to be timelike [resp. lightlike, etc.] if its tangent vectors $\gamma '(t)$ are timelike [resp. lightlike, etc.], for all $t \in I$. Geodesics and geodesic completeness are defined just as the Riemannian case, but geodesics in a Lorentz manifolds are either timelike, lightlike or spacelike, and geodesic incompleteness for each of these types is logically independent of incompleteness for the other types of geodesics \cite{oneill,BE}. Thus, for example, there are Lorentzian manifolds for which all spacelike and lighlike geodesics are complete but that have incomplete timelike geodesics. Finally, a vector field $X: M \rightarrow TM$ is timelike [resp. lightlike, etc.] if $X(p) \in T_pM$ is timelike [resp. lightlike, etc.] for every $p \in M$.  

We say that two {\em timelike}, continuous vector fields $X,Y: M \rightarrow TM$ have {\em the same temporal orientation} if $g(X,Y) <0$. This is an equivalence relation on the set of continuous timelike vector fields on $M$, and it can be shown, using the connectedness of $M$, that there are at most two corresponding equivalence classes. When there are {\em exactly} two such equivalence classes, $M$ is said to be {\em time-orientable}; in this case, $M$ is {\em time-oriented} if one such equivalence class has been chosen. A time-oriented connected Lorentz manifold is called a {\em spacetime}. In what follows we shall always assume that $M$ is a spacetime. In practice, we choose a time-orientation by fixing from now on some fiducial continuous timelike vector field, which we denote by $X_0:M \rightarrow TM$. The simplest example is {\em Minkowski spacetime}, where $M=\mathbb{R}^{n+1}$, and $g$ is given by the following line element in standard coordinates $(x_0, \ldots, x_n)$

\be
\label{minkowski}
ds^2 = -dx_0^2 + \sum_{i=1}^n dx_i^2,
\ee
and the time-orientation is given by $X_0 = \frac{\partial}{\partial x_0}$. 

Let $p \in M$. One can show (see, e.g., pg. 141 of \cite{oneill}) that given $v \in T_pM$ with $g_p (v, X_0(p))=0$, then either $v =0$ or $v$ is spacelike. Hence, if $v$ is nonspacelike, then either $g_p (v, X_0(p))<0$, in which case we shall say that $v$ is {\em future-directed} , or $g_p (v, X_0(p))>0$, and we then say that $v$ is {\em past-directed} \footnote{As the reader will have noticed, most terms used in Lorentzian Geometry have their origin in Physics, and do have natural physical interpretations, but these are not relevant here.}. Accordingly, a differentiable, nonspacelike curve $\gamma:I \subseteq \mathbb{R}\rightarrow M$ is future-directed [resp. past-directed] if $\gamma '(t)$ is future-directed [resp. past-directed], for all $t \in I$. Note that any $C^1$ nonspacelike curve $\gamma:I \subseteq \mathbb{R}\rightarrow M$ must be either future- or past-directed. Given $p,q \in M$, we write $p <<q$ if there exists a smooth future-directed timelike curve segment $\alpha:[0,1] \rightarrow M$ with $\alpha(0) =p, \alpha(1)=q$. Let $A \subseteq M$ be any subset. The {\em chronological future} of $A$ is the set 
\[
I^{+}(A) = \{ q \in M \; : \; p << q \mbox{ for some } p \in A \}.
\]
(The {\em chronological past} is defined time-dually.)   

An embedded, codimension one smooth submanifold $\Sigma \subset M$ is a {\em spacelike hypersurface} if $\Sigma$ with the induced metric is a Riemannian manifold. We shall fix for the rest of this Section one such $\Sigma \subset M$. For technical reasons, we shall assume in addition that every inextendible causal curve intersects $\Sigma$ exactly once. In this case $\Sigma$ is said to be a {\em Cauchy hypersurface} for $M$. (The existence of such a Cauchy hypersurface puts strong constraints in the geometric structure of $M$, and it implies in particular that $\Sigma$ is connected and that $M$ is diffeomorphic to $\Sigma \times \mathbb{R}$ \cite{sanchez}. But it is quite adequate for our purposes here.)

Now, consider a smooth, connected, compact (without boundary), spacelike submanifold $S \subset M$ of codimension two. Suppose $S$ {\em separates} $\Sigma$, i.e., $S \subset \Sigma$ and $\Sigma \setminus S$ is not connected. This means, in particular, that $\Sigma \setminus S$ is a disjoint union $\Sigma_{+}\dot{\cup} \Sigma_{-}$ of open submanifolds of $\Sigma$ having $S$ as a common boundary. We shall loosely call $\Sigma_{+}$ [resp. $\Sigma_{-}$] the {\em outside} [resp. {\em inside}] of $S$ in $\Sigma$. (In most interesting examples there is a natural choice for these.) It also means that $S$ is two-sided in $\Sigma$ and thus there are unique unit spacelike vector fields $Z_{\pm}$ on $S$ normal to $S$ in $\Sigma$, such that $Z_{+}$ [resp. $Z_{-}$] is outward-pointing [resp. inward-pointing]. 

Note that because $M$ is time-oriented, there exists a unique timelike, future-directed, normal vector field $U$ on $\Sigma$ satisfying $g(U,U) =-1$. Thus the normal bundle of $\Sigma$ is trivial, and $K_{\pm} := U|_{S}+Z_{\pm}$ are future-directed {\em lightlike} vector fields on $S$ normal to $S$ in $M$. The {\em outward }[resp. {\em inward}] {\em null expansion scalar} of $S$ in $M$ is the smooth function $\theta_{+}: S \rightarrow \mathbb{R}$ [resp. $\theta_{-}: S \rightarrow \mathbb{R}$] given by 
\begin{equation}
\label{maineq1}
\theta _{+}(p) = - g_p( H_p, K_{+}(p)) \mbox{  [resp. $\theta_{-}(p) = - g_p (H_p, K_{-}(p)) $]},
\end{equation}
for each $p \in S$, where $H_p$ denotes the mean curvature vector of $S$ in $M$ at $p$ \cite{oneill}. 

Geometrically, these functions measure the initial divergence of families of future-directed lightlike normal geodesics emanating from $S$. (Since $S$ has codimension 2 there are exactly two such families, whose initial tangent vectors at a point $p \in S$ are parallel to $K_{\pm}(p)$.) If $S$ is a round sphere in a Euclidean slice $x_0 = \mbox{ const.}$ of Minkowski spacetime (\ref{minkowski}), with the obvious choices of inside and outside, then we have $\theta_{+}>0$ and $\theta_{-}<0$. 

Using the terminology above, we end this Section with the following useful definition:
 
\begin{definition}
\label{asymptoticallyregular}
 The spacelike Cauchy hypersurface $\Sigma \subset M$ is {\em asymptotically regular} if there exists a smooth, connected, compact embedded submanifold $S \subset \Sigma$ of dimension $n-1$ such that
\begin{itemize}
\item[i)] $S$ separates $\Sigma$, and $\overline{{\Sigma}}_{+} \equiv S \cup \Sigma_{+}$ is noncompact;
\item[ii)] The map $j_{\#}: \pi_1(S) \rightarrow \pi_1(\overline{\Sigma}_{+})$ induced by the inclusion $j: S \hookrightarrow \Sigma$ is onto; 
\item[iii)] $\theta_{-} <0$ everywhere on $S$.
 \end{itemize}
we shall call such an $S$ an {\em enclosing surface} in $\Sigma$. 
\end{definition}

%===========================================================
\section{An Application to Asymptotically Flat Manifolds}\label{AsympFlat}
%==========================================================

Throughout this Section, $N=(N,h)$ is a complete, noncompact, connected $n$-dimensional ($n \geq 2$) Riemannian manifold with nonnegative Ricci curvature $Ric_N\geq 0$. 

The main motivation of this Section is to see that the assumptions of our main theorem are naturally met in some concrete settings. One especially important such setting is when $N$ is {\em asymptotically flat}. This concept arises quite naturally in General Relativity, specifically in the study of gravitationally isolated systems (see, e.g., \cite{wald} for a introduction to the physical aspects of the subject); but it has attracted the continued attention of geometers at least since Schoen-Yau's celebrated proof of the Positive Mass Theorem \cite{SY1,SY2}. This theorem establishes that the so-called Arnowitt-Deser-Misner (ADM) mass of asymptotically flat manifolds of nonnegative scalar curvature is a nonnegative number, and provides a powerful global geometric invariant: if this mass is zero, then the manifold is {\em isometric} to $E^n$. Schoen-Yau arguments, given initially for $n=3$, work for dimensions $3\leq n \leq 7$; moreover, Witten \cite{witten} has given a proof of the result which works \cite{bartnik} for all dimensions in the particular case when the manifold is spin. 

Let us first define more precisely what we mean here by `asymptotically flat'. 

\begin{definition}
\label{AFdef}
$(N,h)$ is {\em asymptotically flat} if there exists a compact subset $\Omega \subset N$ such that $M\setminus \Omega$ has a finite number of components $E_1, \ldots, E_k$, called {\em (asymptotically flat) ends}, and each such end $E_i$ is diffeomorphic to the region $\{ x \in \mathbb{R}^n \; : \; |x| >1 \}$. In addition, the metric components in the coordinate system induced by this diffeomorphism satisfy the estimates
\[
|h_{ij}(x) - \delta_{ij}| \leq \frac{C}{|x|^{\alpha}} \mbox{     and     }
|\partial _k h_{ij}(x)| \leq \frac{C}{|x|^{\alpha +1}},
\]
for some positive constants $C,\alpha$. \footnote{The definition of asymptotic flatness varies slightly throughout the literature. The one presented here is weaker than the ones used, e.g., in Refs. \cite{SY1,SY2, bartnik}. We believe, however, that our requirements are minimal in the sense that manifolds which are asymptotically flat in other definitions must also be so according to ours.} 
\end{definition}

If $N$ is asymptotically flat, then since $N$ has been assumed to be complete with nonnegative Ricci curvature, it is a well-known consequence of the Cheeger-Gromoll Splitting theorem \cite{CG} that it must have at most two asymptotically flat ends. We show that it has exactly one end. 

\begin{theorem}
\label{AFtheo}
If $N$ is a complete connected $n$-dimensional ($n \geq 3$) asymptotically flat Riemannian manifold with Ricci curvature $Ric_N\geq 0$, then $N$ is simply connected and has exactly one end. 
\end{theorem}

{\em Proof.} Let $E \subset N$ be an asymptotically flat end, and let 
\[
\psi : R=\{ x \in \mathbb{R}^n \; : \; |x| >1 \} \rightarrow E
\]
be a diffeomorphism. Let $r>1$ and $S_r = \{ x \in \mathbb{R}^n \; : \; |x| = r \} \subset R$. The image of $S_r$ by $\psi$ (which by a slight abuse of notation we also denote as $S_r$ in what follows) is an embedded hypersurface contained in $E$. If we write 
\[
N_{+} = \psi\left(\{ x \in \mathbb{R}^n \; : \; |x| >r \}\right), \\
N_{-} = N \setminus \psi\left(\{ x \in \mathbb{R}^n \; : \; |x| \geq r \}\right),
\]
then, clearly, $S_r \cup N_{+}$ is noncompact and $i_{\#}:\pi_1(S_r) \rightarrow \pi_1(S_r \cup N_{+})$ induced by the inclusion is surjective. We now proceed show that {\em for some choice of $r>0$, $S_r$ has positive mean curvature} (with respect to the unit normal pointing to $N_{+}$). In this case, by Theorem \ref{maintheorem}, $S_r \cup N_{-}$ is compact, which means that $N$ can have no other asymptotically flat end, and $i_{\#}:\pi_1(S_r) \rightarrow \pi_1(N)$ is surjective, and hence $\pi_1(N) =0$ since $n\geq 3$. 

Thus, let $f_0:R \rightarrow \mathbb{R}$ be given by $f_0(x) = |x|$, and $f:= f_0 \circ \psi^{-1}$. Let $h_0$ be the flat metric induced on $E$ by $\psi^{-1}$, and let $\nabla f$, $\nabla _0 f$ be the gradients of $f$ with respect to $h$ and $h_0$ respectively. Of course, $S_r$ is a level set for $f$. 

Fix $p \in S_r$. Denote by $\langle \; , \; \rangle$ and $\langle \; , \; \rangle _0$ the inner products on $T_pS_r$ defined by $h$ and $h_0$ respectively, and let $B: T_pS_r \rightarrow T_pS_r$ be the unique self-adjoint (with respect to $\langle \; , \; \rangle _0$) operator such that 
\[
\langle v , w \rangle = \langle v , Bw \rangle _0 , \forall v,w \in T_pS_r.
\]
Let $\{e_1, \ldots, e_{n-1} \}$ be an orthonormal basis (with respect to $\langle \; , \; \rangle _0$) of $B$ eigenvectors, and write $B e_i = \lambda_i e_i$ for the eigenvalues; of course $\lambda_i >0$ ($i=1, \ldots, n-1$). Therefore $\{ \frac{e_1}{\sqrt{\lambda_1}}, \ldots, \frac{e_{n-1}}{\sqrt{\lambda_{n-1}}} \}$ is orthonormal with respect to $\langle \; , \; \rangle$. Finally, choose $V_1 , \ldots, V_{n-1}$ smooth vector fields on $E$ such that $V_i(p) = \frac{e_i}{\sqrt{\lambda_i}}$ ($i=1, \ldots, n-1$). For the remaining part of the proof the computations are understood to be carried out at $p$. 

The mean curvature $H = H_{S_r}$ at $p$ is given by
\[
H = - \frac{1}{n-1} \sum_{i=1}^{n-1} h \left( \nabla ^h_{V_i}V_i, \frac{\nabla f}{| \nabla f|} \right) \equiv - \frac{1}{n-1} \sum_{i=1}^{n-1} \frac{1}{| \nabla f|}h_0 \left( \nabla ^h_{V_i}V_i, \nabla _0 f \right), 
\] 
where $|\nabla f| = \sqrt{\langle \nabla f (p), \nabla f(p) \rangle }$. 

Introduce, for any two smooth vector fields $X,Y$ on $E$, the tensorial quantity 
\[
D(X,Y) := \nabla^h_X Y - \nabla^{h_0}_X Y. 
\]
Using the coordinate system on $E$ defined by $\psi$ we have
\[
D(\sum_{i=1}^{n} X^i \partial_i , \sum_{j=1}^{n} Y^j \partial_j ) = \frac{1}{2} \sum_{i,j,k,l =1}^n X^iY^j h^{kl}\left( \partial_j h_{il} + \partial_i h_{lj} - \partial_l h_{ij}\right) \partial_k ,
\]
where here and hereafter we denote by $h^{ij}$ the components of the inverse matrix of $[h_{ij}]$. If we write $e_i = \sum_{k=1}^{n} e^k_i \partial_k |_p$, a straightforward computation gives that
\be
\label{mean}
H = \frac{1}{(n-1)|\nabla f|} \sum_{i=1}^{n-1} \left( \frac{1}{ \lambda_i r} - \langle D(V_i,V_i), \nabla_0 f \rangle_0 \right),
\ee
where 
\[
\langle D(V_i,V_i), \nabla_0 f \rangle_0 = \sum_{i=1}^{n-1} \sum_{k,l,a,b =1}^n \frac{1}{2} \frac{e_i^k e_i^l x^a}{r \lambda_i}h^{ab}\left( \partial_l h_{bk} + \partial_k h_{bl} - \partial_b h_{kl} \right). 
\]
Moreover it is easy to check that $|e^k_i| \leq 1$ ($1 \leq i \leq n-1, 1 \leq k \leq n$). 

Now we write 
\[
h_{ij} = \delta_{ij} + t_{ij} \; , \; h^{ij} = \delta_{ij} + t^{ij}, 
\]
whence we have $t^{ij} = - \sum_k t_{jk} t^{ik} - t_{ij}$. Recall that for some numbers $C, \alpha >0$, 
\[
|t_{ij}(x)| \leq \frac{C}{|x|^{\alpha}} \mbox{     and     }
|\partial _k h_{ij}(x)| \leq \frac{C}{|x|^{\alpha +1}}
\]
for $|x| >1$. We now assume for later convenience that 
\begin{equation}
\label{estimate}
r^{\alpha}  > 6n^4C .
\end{equation}
From (\ref{estimate}) we easily check that, since in particular $r^{\alpha}  > 2nC$, 
\[
|t^{ij}| \leq  \sum_{k=1}^n |t_{jk}|| t^{ik}| + |t_{ij}| < \frac{2C}{r^{\alpha}}.
\]
Thus for $i=1, \ldots, n-1$,
\begin{eqnarray}
|\langle D(V_i,V_i), \nabla_0 f \rangle_0 | &\leq & \frac{1}{2 \lambda_i r} \sum_{k,l,a,b =1}^n |e^k_i||e^l_i||x^a|(\leq r)|h^{ab}|\left( |\partial_ l h_{bk}| + |\partial_k h_{bl}| + |\partial _b h_{kl}| \right) \nonumber \\ 
&\leq & \frac{1}{2\lambda_i r}\left( 1 + \frac{2C}{r^{\alpha}}\right)\left( \frac{3Cn^4}{r^{\alpha}} \right) < \frac{1}{2\lambda_i r}\left( \frac{6Cn^4}{r^{\alpha}} \right) <\frac{1}{\lambda_i r}, \nonumber
\end{eqnarray}
again because of the estimate (\ref{estimate}). Thus, using Eq. (\ref{mean}) we conclude that $H>0$, as desired. 
\qcd

{\em Remark 1.} The assumption of completeness in Theo. \ref{AFtheo} cannot be dropped. Indeed, if completeness fails, the manifolds $\mathbb{R}^3 \setminus \{ (x,y,0) \in \mathbb{R}^3 \; : \; x^2 + y^2 = 1 \}$ and $\mathbb{R}^3 \setminus \{ (x,y,z) \in \mathbb{R}^3 \; : \; x^2 + y^2 + z^3 \leq 1 \}$, both with a flat metric, are counterexamples for the the simple-connectedness conclusion and for the existence of one single end, respectively. 

{\em Remark 2.} The assumption that $Ric_N \geq 0$ in Theo. \ref{AFtheo} cannot be weakened to nonnegative {\em scalar} curvature $R_N \geq 0$. As a simple example, consider the {\em Schwarzschild 3-manifold of mass $m$} ($m \in (0, +\infty)$) $N_S= \mathbb{R}^3  \setminus 0$ with the metric given by 
\[
h_{ij}(x)= \left( 1 + \frac{m}{2 |x|} \right)^4 \delta_{ij}
\]
for all $x \in N$, where $m$ is some positive real constant. One can check by direct calculation that $N_S$ is complete and has nonnegative scalar curvature, but its Ricci tensor is indefinite. Now, $N_S$ has two asymptotically flat ends. One flat end clearly corresponds to the region of large $|x|$. But the region around the ``origin'' $x=0$ also corresponds to an asymptotically flat end. To see this explicitly, note that $N_S$ is isometric to $(0, +\infty) \times S^2$ with the metric given by the line element
\[
ds^2 = \left( 1 + \frac{m}{2 \rho} \right)^4 \left( d\rho ^2 + \rho^2 d\Omega^2\right), 
\]
where $d\Omega ^2$ is the standard metric on $S^2 = \{p \in \mathbb{R}^3 \; : \; |p| =1 \}$. The transformation $(\rho, p) \mapsto (\frac{m^2}{4\rho},p)$ for all $\rho \in (0, +\infty)$ and all $p \in S^2$ is an isometry which swaps the two asymptotically flat regions. Moreover, the mapping $(\rho, p) \mapsto (\frac{m^2}{4\rho},-p)$ for all $\rho \in (0, +\infty)$ and all $p \in S^2$ is also an idempotent isometry which generates a $\mathbb{Z}_2$-subgroup of the isometry group of $N_S$ which acts properly discontinuously on $N_S$. Quotienting $N_S$ by this action gives a manifold $\tilde{N}_S$ locally isometric to $N_S$ (and so it also has nonnegative scalar curvature and indefinite Ricci tensor), but which is diffeomorphic to the projective space $\mathbb{R}P^3$ minus a point, and thus it is not simply connected.  

{\em Remark 3.} Given additional geometric requirements, the fact that a complete manifold is simply connected can be instrumental in proving more stringent results. For example, Nardmann \cite{nardmann} proves that $N$ can be embedded in either Minkowski spacetime or in (the universal covering of) anti-de Sitter spacetime as a graph (which means that $N$ is diffeomorphic to $\mathbb{R}^n$) if certain intrinsic analogues of the Gauss-Codazzi equations are satisfied on $N$, together with some additional topological requirements. These latter requirements are automatically satisfied if $N$ is a complete connected $n$-dimensional ($n \geq 3$) asymptotically flat Riemannian manifold with Ricci curvature $Ric_N\geq 0$, for in this case $N$ is simply connected in view of our result. Thus, if in addition the mentioned Gauss-Codazzi analogues are satisfied in this case, the conclusions of the main theorem of Ref. \cite{nardmann} hold. We refer the reader to that reference for further details. 

{\em Remark 4.} For $n=3$, the hypotheses of the Theo. \ref{AFtheo} imply that $N$ is actually diffeomorphic to $\mathbb{R}^3$. In 1994, Zhu \cite{zhu}, improving on a previous result by Schoen and Yau \cite{SY3}, showed that a complete $3$-manifold of nonnegative Ricci curvature (not necessarily asymptotically flat) is diffeomorphic to $\mathbb{R}^3$ if $Ric_N >0$ at least at one point. But if this fails, then $N$ is flat. In that case, if in addition $N$ minus some compact set is diffeomorphic to $\mathbb{R}^3$ minus a ball, then the associated ADM mass is zero and $N$ is diffeomorphic to $\mathbb{R}^3$ by the Positive Mass theorem. Theo. \ref{AFtheo} offers another proof of this fact.

%===========================================================
\section{Proof of the Main Theorem}\label{detailedproofs}
%==========================================================

In a previous paper by this author \cite{GLme}, the following two results, which are given here as Lemmas \ref{auxiliary} and \ref{maintheoremprevious}, have been proven \footnote{The notation in the mentioned paper is somewhat different, and the hypotheses used there are actually weaker than those stated here, but the present versions suffice for our purposes.}. Since the full proofs are given in that reference, only sketches are presented here. Throughout this Section, $(N,h)$ is a complete, noncompact, connected $n$-dimensional ($n \geq 2$) Riemannian manifold with Ricci curvature $Ric_N\geq 0$. 

\begin{lemma}
\label{auxiliary}
Let $(M,g)$ be an $(n+1)$-dimensional (with $n \geq 2$) lightlike geodesically complete spacetime, which satisfies $Ric_M(v,v) \geq 0$ for all lightlike $v \in TM$ and admits an asymptotically regular Cauchy hypersurface $\Sigma \subset M$. Then, given an enclosing surface $S \subset \Sigma$, the closure of the inside of $S$, $\overline{\Sigma}_{-} \equiv S \cup \Sigma_{-}$, is compact. 
\end{lemma}

{\em Sketch of the Proof.} The idea is to show, first, that the set $T:= \partial I^{+}\left( \Sigma_{+} \right) \setminus \left( \Sigma_{+} \right)$ is compact. This follows by noting that $T$ is generated by the family of future-directed normal null geodesics starting at $S$ with initial tangent vector parallel to $K_{-}$: since $\theta_{-} <0$ everywhere on $S$, $Ric_M(v,v) \geq 0$ for all lightlike $v \in TM$ and $M$ is null geodesically complete, standard considerations in the causal theory of spacetimes imply that $T$ must be compact (see, e.g., pg. 436 of \cite{oneill}).  Now, since $\Sigma$ is a spacelike Cauchy hypersurface, all inextendible integral curves of the fiducial timelike vector field $X_0$ which intersect $T$ must also intersect $\Sigma$, and one can check that this induces a homomorphism of $T$ onto $\overline{\Sigma}_{-}$. 
\qcd

\begin{lemma}
\label{maintheoremprevious}
Let $(M,g)$ be an $(n+1)$-dimensional (with $n \geq 2$) lightlike geodesically complete spacetime, which satisfies $Ric_M(v,v) \geq 0$ for all lightlike $v \in TM$ and admits an asymptotically regular Cauchy hypersurface $\Sigma \subset M$. Then, given an enclosing surface $S \subset \Sigma$, the group homomorphism $j_{\#}: \pi_1(S) \rightarrow \pi_1(\Sigma)$ induced by the inclusion $j: S \hookrightarrow \Sigma$ is surjective. In particular, if $S$ is simply connected, then so is $\Sigma$.  
\end{lemma}

{\em Sketch of the Proof.} The strategy here is to consider a certain smooth connected covering $\phi: \tilde{M} \rightarrow M$ endowed with the pullback metric $\tilde{g} := \phi^{\ast} g$ and induced time orientation so that it is a spacetime by itself. The main features of this specific covering are: (i) $(\tilde{M},\tilde{g})$ is an $(n+1)$-dimensional lightlike geodesically complete spacetime, which satisfies $Ric_{\tilde{M}}(v,v) \geq 0$ for all lightlike $v \in T\tilde{M}$ and admits an asymptotically regular Cauchy hypersurface $\tilde{\Sigma} \subset \tilde{M}$, and (ii) there exists an enclosing surface $\tilde{S} \subset \tilde{\Sigma}$ diffeomorphic to $S$ and with $\tilde{\Sigma}_{+} \cup \tilde{S}$ diffeomorphic to $\Sigma_{+} \cup S$. Such a covering always exists \cite{hawcovering} and is called sometimes {\em the Hawking covering} in the literature \footnote{I am grateful to Gregory Galloway for pointing this out to me, as well as for bringing Ref. \cite{hawcovering} to my attention.}. If $j_{\#}: \pi_1(S) \rightarrow \pi_1(\Sigma)$ is not onto, then a detailed analysis reveals that $\tilde{S} \cup \tilde{\Sigma}_{-}$ is noncompact, contradicting the previous Lemma. 
\qcd

The latter result is a generalization of some classic theorems obtained independently by Gannon \cite{gannon1,gannon2} and Lee \cite{lee} in the 1970's (see \cite{GLme} for a detailed discussion).  

\vspace{0.5cm}

{\em Proof of Theorem \ref{maintheorem}}:
\newline
Consider the following class of spacetimes: we take $M = \mathbb{R} \times N$ with the Lorentz metric $g_f := -dt^2 \oplus f^2 h$, where $f:\mathbb{R} \rightarrow \mathbb{R}$ is some smooth positive function, time-oriented by the timelike vector field $\frac{\partial}{\partial t}$. This is a particular case of a Lorentzian warped product \cite{oneill,BE}. For a fixed $t_0 \in \mathbb{R}$, the smooth map $\varphi_{t_0}: N \rightarrow M$ given by $\varphi_{t_0}(p)= (t_0,p)$ for each $p \in N$ is an embedding (it is {\em not} in general an isometric embedding; rather  $\varphi_{t_0}^{\ast} g_f = f^2(t_0) h$), and therefore the image $\Sigma_{t_0}= \varphi_{t_0}(N)$ of $N$ by $\varphi_{t_0}$ can be viewed as a spacelike hypersurface \cite{oneill,BE} in $M$. 

Now, the hypotheses of Theorem \ref{maintheorem} just mean that $N$ admits an embedded hypersurface $S \subset N$ which satisfies conditions $(i)$ and $(ii)$ of Definition \ref{asymptoticallyregular}. However, a computation of the null expansion scalars $\theta^{t_0}_{\pm}$ of $S_{t_0} = \varphi_{t_0}(S)$ in this case gives (see also, e.g., Eq. (2.2) of Ref. \cite{flores})
\[
\theta^{t_0}_{\pm} = \frac{f'(t_0)}{f(t_0)} \pm \frac{H_S}{f(t_0)},
\]
where $H_S$ is the mean curvature of $S$ in $N$ with respect to the ambient metric $h$ and the normal outside-pointing vector field $Z_{+}$. 

Choose $f\equiv 1$ and write $g_1 \equiv g$. In this case $\theta_{-}^{t_0} = -H_S <0$, since we assume $H_S >0$, and therefore $\Sigma_{t_0}$ can be viewed in this case as an asymptotically regular {\em Cauchy} hypersurface in $M$. (That is because $(N,h)$ is complete. Moreover, it is well-known that $(M,g)$ is geodesically complete in this case - see, e.g., Theo. 3.67, pg. 103, of Ref. \cite{BE} for a proof of these facts.) Now, given a future-directed lightlike vector $v \in TM$, it can be easily checked that it must be of the form $v= \lambda \left(\frac{\partial}{\partial t} + u_{L} \right)$, where $\lambda >0$ and $u_{L}$ is the lift to $M$ of a unit vector $u \in TN$. Standard formulas for the Ricci tensor of a warped product (see, e.g., pag. 211 of \cite{oneill}) give that 
\[
Ric_M( \frac{\partial}{\partial t},\frac{\partial}{\partial t})= Ric_M(\frac{\partial}{\partial t},u_L) =0 \mbox{ and } Ric_M(u_L,u_L) = Ric_N(u,u) \geq 0,
\]
since the Ricci tensor of $(N,h)$ is assumed to be nonnegative. Therefore $Ric_M(v,v) \geq 0$, and $(M,g)$ satisfies all the conditions of Lemmas \ref{auxiliary} and \ref{maintheoremprevious}. We therefore conclude that $S \cup N_{-}$ is compact and that $i_{\#}:\pi_1(S) \rightarrow \pi_1(N)$ is surjective, so the proof is complete.
\qcd

\newpage  
%%%%%%%%%%%%%%%%%%%%%%%%%%%%%%%%%%%%%%%%%%%%%%%%%%%%%%%%%%%%%%%%%%%%%%%%%%%%%%%%%%%%% 
%References 
%%%%%%%%%%%%%%%%%%%%%%%%%%%%%%%%%%%%%%%%%%%%%%%%%%%%%%%%%%%%%%%%%%%%%%%%%%%%%%%%%%%%%

%%%%%%%%%%%%%%%%%%%%%%%%%%%%%%%%%%%%%%%%%%%%%%%%%%%%%%%%%%%%%
%
%%%%%%%%%%%%%%%%%
\end{document}